\documentclass[a4paper]{amsart}
\usepackage{amsmath,amsfonts,amsthm,amssymb,color,bm,hyperref}
\numberwithin{equation}{section}
\newtheorem{theorem}{Theorem}[]

\newtheorem{lemma}[theorem]{Lemma}

\newcommand*{\abs}[1]{\lvert#1\rvert}

\begin{document}

\title[Normal subgroups and permutation characters]{Normal subgroups and permutation characters: a correction to a proof of Klingen}  

\author[P.~M\"{u}ller]{Peter M\"{u}ller}
\address{
Institut f\"ur Mathematik, Universit\"at W\"urzburg, 97074 Würzburg, Germany
}
\email{peter.mueller@uni-wuerzburg.de}

\author[P.~Spiga]{Pablo Spiga}
\address{Dipartimento di Matematica e Applicazioni, University of Milano-Bicocca, Via Cozzi 55, 20125 Milano, Italy} 
\email{pablo.spiga@unimib.it}
\begin{abstract}
  Let $G$ be a finite group. For subgroups $U$ and $V$ let $1_U^G$ and
  $1_V^G$ be the permutation characters for the action of $G$ on the
  right cosets of $U$ and $V$, respectively. Let $N$ be a normal
  subgroup of $G$.  Norbert Klingen, in his book \cite{klingen}, shows
  that if $1_U^G=1_V^G$, then $1_{UN}^G=1_{VN}^G$. We give a
  counterexample to an argument in his proof and provide a correct
  proof for this statement.
\end{abstract}


\maketitle

\section{Introduction}\label{sec:intro}
Many questions on Kronecker classes, and more generally in algebraic
number theory, can be phrased in purely group-theoretic
terms. Examples can be found, for instance, in the work of
Jehne~\cite{Jehne77} and Perlis~\cite{perlis}. Norbert Klingen, in his beautiful book~\cite{klingen}, further explores these connections and provides a comprehensive treatment of this area, offering a monograph that remains highly valuable today.

One particularly relevant situation that arises in the study of arithmetic similarities is the equality $1_U^G = 1_V^G$ between permutation characters. In this context, Klingen states the following result.

\begin{theorem}[Theorem~1.6~(a), page 83,~\cite{klingen}]\label{thrm}Let $G$ be a finite group, let $U$ and $V$ be  subgroups of $G$, 
let $1_U^G$ and $1_V^G$ be the permutation characters for the action of $G$ on the right cosets of $U$ and $V$ respectively and  let $N$ be a normal subgroup. If $1_U^G=1_V^G$, then $1_{UN}^G=1_{VN}^G$.
\end{theorem}
However, the proof of this theorem is not correct. Indeed, the smallest counterexample to the essential step in Klingen's proof is the following (we use the notation from~\cite[page~83]{klingen}, except that we prefer to use the letter $N$ to denote the normal subgroup $H$ in the proof of ~\cite[Theorem~1.6~(a)]{klingen}). Let $G = \langle x \rangle$ be cyclic of order $4$, let $U$ be the identity subgroup of $G$, let $N = \langle x^2 \rangle$, and let $\sigma = x^2$. Then $\sigma \in UN$, so $1_{UN}^G(\sigma) > 0$, while $1_U^G(\sigma) = 0$. Therefore,~\cite[equality~(14)]{klingen} does not hold.

Theorem \ref{thrm}, which is an immediate consequence of the following
Lemma \ref{lemm}, is very relevant for inductive arguments.
\begin{lemma}\label{lemm}
  Let $U$ be a subgroup of the finite group $G$, and let $N$ be a
  normal subgroup of $G$. Then
  \[
    1_{UN}^G(g)=\frac{1}{\abs{N}}\sum_{n\in N}1_U^G(gn)\text{ for all
    }g\in G.
    \]
\end{lemma}
\section{Proof of Lemma \ref{lemm}}
Lemma \ref{lemm} is a version of Lemma \ref{fgs} below which first
appeared in the work by Fried, Guralnick and Saxl on exceptional
permutation polynomials, see \cite[proof of Lemma 13.1]{FGS}. Its
original proof employed representation theoretic arguments. It was
later proved again in \cite[Lemma 3.1]{GW} by a short counting
argument.
\begin{lemma}[Fried, Guralnick, Saxl]\label{fgs} Let the finite group
  $H$ act on the finite set $\Omega$, let $N$ be a normal subgroup of $H$,
  let $\pi(h)$ be the number of fixed points of $h\in H$ on $\Omega$, and
  assume that $H=\langle N,g\rangle$ for some $g\in H$. Let $r$ be the
  number of $H$-orbits on $\Omega$ which do not split into smaller
  $N$-orbits. Then
  \[
    \frac{1}{\abs{N}}\sum_{n\in N}\pi(gn)=r.
  \]
\end{lemma}  
Note that Lemma \ref{lemm} follows from this result: Let $\Omega$ be
set of right cosets of $U$ in $G$ and set $H=\langle N,g\rangle$. Then
$1_U^G=\pi$ and $1_{UN}^G(g)$ is the number of $N$-orbits on $\Omega$
fixed by $g$. As $H=\langle N,g\rangle$, these are exactly the
$H$-orbits which do not split into smaller $N$-orbits.
\section*{Acknowledgments}
 The second author is funded by the European Union via the Next
Generation EU (Mission 4 Component 1 CUP B53D23009410006, PRIN 2022, 2022PSTWLB, Group
Theory and Applications).
\end{document}